\documentclass[leqno,11pt]{amsart}
\usepackage{amsmath}
\usepackage{amssymb,amscd}
\usepackage{amsthm}
\usepackage{latexsym}
\usepackage[myheadings]{fullpage}
\usepackage{color}  

\newtheorem{theorem}{Theorem}[section]
\newtheorem{lemma}[theorem]{Lemma}

\newtheorem{proposition}[theorem]{Proposition}

\newtheorem{remark}[theorem]{Remark}
\newtheorem{example}[theorem]{Example}

\newcommand{\m}{\mathfrak m}
\newcommand{\n}{\mathfrak n}

\newcommand{\Hom}{{\rm Hom}}

\def\opn#1#2{\def#1{\operatorname{#2}}} 
\opn\spec{Spec}
\opn\depth{depth}
\opn\height{ht}
\opn\chara{char}
\opn\gr{gr}
\opn\ord{ord}
\opn\Ap{Ap}
\opn\F{F}
\opn\PF{PF}
\opn\L{L}
\opn\Max{Max}

\title{Some classes of one-dimensional rings characterized by their reflexive ideals}

\author{Pietro Campochiaro}
\address{Pietro Campochiaro - Scuola Superiore di Catania - Universit\`a degli Studi di Catania - Via Val\-disavoia 9 - 95123 Catania - Italy}
\email{campochiaropietro.18@gmail.com}

\author{Marco D'Anna}
\address{Marco D'Anna - Dipartimento di Matematica e Informatica - Universit\`a degli Studi di Catania - Viale Andrea Doria 6 - 95125 - Catania - Italy}\email{mdanna@dmi.unict.it}

\author{Francesco Strazzanti}
\address{Francesco Strazzanti - Dipartimento di Scienze Matematiche e Informatiche, Scienze Fisiche e Scienze della Terra - Universit\`a degli Studi di Messina - Viale Ferdinando Stagno d’Alcontres 31 - 98166 - Messina - Italy}
\email{francesco.strazzanti@unime.it}

\thanks{
The second author was supported by the PRIN 2020 ``Squarefree Gröbner degenerations, special varieties and related topics''. 
The second and third authors are members of the ``National Group for Algebraic and Geometric Structures, and their Applications'' (GNSAGA - INdAM)}


\date{}

\begin{document}

\begin{abstract}
We study reflexive ideals in one-dimensional Cohen-Macaulay local rings, providing characterizations of almost Gorenstein rings, rings with minimal multiplicity, and Arf rings, which describe their reflexive fractional ideals.  
\end{abstract}


 \maketitle

\section{Introduction}
Let $R$ be a commutative ring with unit element. An $R$-module $M$ is said to be \emph{reflexive} if $M\cong \Hom_R(\Hom_R(M,R),R)$.
When $I$ is a regular fractional ideal of $R$, we have
$\Hom_R(I,R)\cong (R: I)$, where the colon is taken in the total ring of fractions $Q(R)$ of $R$, and therefore reflexivity can be expressed as $I=R:(R:I)$.
In this case, $I$ is also said \emph{divisorial}.

Reflexivity of modules and ideals and its interactions with the properties of the ring is a classical subject of studies, especially in the context of integral domains (see e.g. \cite{Jans} and \cite{M0}).

Recently, many papers have focused on the reflexivity problems (about both ideals and modules) for one-dimensional, noetherian rings coming from curve singularities (see, for instance, \cite{D}, \cite{DMS}, \cite{IK}, \cite{K}).
In \cite{DMS}, the authors point out why the one-dimensional case is particularly significant. In fact, for noetherian rings, reflexivity of a finitely generated module $M$ over a noetherian $(S_2)$ ring is equivalent to $M$ satisfying $(S_2)$ and being reflexive in codimension one, see \cite[Proposition 1.4.1]{BH}.

In the same paper \cite{DMS}, the authors address several questions about the nature and the number of reflexive ideals and modules, and, in particular, provide
sufficient conditions for reflexivity of ideals, introducing the notion of $I$-Ulrich module. It turns out that the 
conductor $R:\overline R$ (where $\overline R$ is the integral closure of $R$ in $Q(R)$), the existence of canonical ideals, and birational extensions play a key role in the study of reflexive ideals and modules. 

In this paper, we focus on (fractional) ideals and aim to deepen the understanding of reflexive ideals, working within a slightly more restrictive setup than that of \cite{DMS}. In fact, we will assume that $R$ is a one-dimensional, Cohen-Macaulay local ring with infinite residue field, and we will assume that $R$ is analytically unramified.
In this way, we are sure that $R:\overline R$ is a regular ideal of $R$, and that every ideal has a principal reduction. These two facts will play a fundamental role in our results. We will also make use of canonical ideals, that exist since $\hat R$ is reduced. In Section 3, using these ingredients, we are able to restrict our attention to ideals $I$ lying between $R:\overline R$ and $R$. For these ideals, we compare divisorial closure to integral closure
(see Proposition \ref{main reflexivity}), obtaining a new necessary condition for the reflexivity of an ideal $I$, which also involves the birational extension $I:I$ (see Lemma \ref{I:I-2}).

Changing point of view, in the famous paper \cite{B},  
Bass proved, among many other results, that a noetherian ring $R$ with $\dim(R)\leq 1$ is Gorenstein if and only if it is Cohen-Macaulay and every ideal $I$ of $R$ is reflexive.
This result shows that the knowledge of reflexive ideals gives information about the ring, or, conversely, that for some classes of rings it is possible to determine their reflexive ideals. Thus a natural question arises:
are there other relevant classes of rings that can be characterized looking at their reflexive ideals?
Trying to answer this question, we consider three important classes: 
almost Gorenstein rings, rings with minimal multiplicity and Arf rings.

Almost Gorenstein rings were firstly introduced in \cite{BF} in the one-dimensional analytically unramified case, and the definition was successively extended to general one-dimensional Cohen-Macaulay rings in \cite{GMP}, with the aim of finding a class of Cohen-Macaulay rings that are close to be Gorenstein under certain respects. This notion has been widely explored in recent years also in higher dimension (see e.g.  \cite{GTT}) and revealed to be very interesting from several points of view. Assuming our setup and denoting by $\omega_R$ a canonical ideal of $R$ with the property that $R\subseteq \omega_R\subseteq \overline R$, we can define almost Gorenstein local rings via the inclusion $\omega_R\m\subseteq \m$
(where $\m$ is the maximal ideal of $R$). In Section 4, we give two characterizations of almost Gorenstein rings involving reflexive ideals and either $\omega_R$ (cf. Theorem \ref{almost Gorenstein}) or 
the birational extension $\m:\m$
(cf. Theorem \ref{almost Gorenstein 2}). A key fact in the second characterization is another necessary and sufficient condition for the reflexivity of an ideal $I$ lying between $R:\overline R$ and $R$, which involves integrally closed ideals containing $I$ (cf. Proposition \ref{H-reflexivity}).

Finally, in Section 5, we study rings with minimal multiplicity and Arf rings. It is well known that, for Cohen-Macaulay local rings of dimension $d$, it holds true that $\nu(R) \leq e(R)+d-1$, where $e(R)$ denotes the multiplicity of $R$ and $\nu(R)$ its embedding dimension, see \cite{A}.
This bound is sharp and rings for which the equality is attained are called rings of minimal multiplicity (or of maximal embedding dimension)
and have particularly nice properties (see e.g. \cite{S}). In our setup, having minimal multiplicity is equivalent to the stability of the maximal ideal $\m$, i.e. $\m^2=x\m$ (where $(x)$ is any minimal reduction of $\m$), or equivalently, $\m:\m=x^{-1}\m$; moreover, if this is the case, the blowing up of $R$ coincides with $\m:\m$. In Theorem \ref{minimal multiplicity}, we provide a characterization of these rings, showing that their reflexive fractional ideals are exactly the reflexive fractional ideals of $\m:\m$. Additionally, $R$ has minimal multiplicity if and only if $\m$ is $(\m:\m)$-reflexive.

Arf rings were defined in \cite{L} in connection with curve singularities; they play a relevant role in the 
study of equisingularity problems, since in the same paper it is proven that the multiplicity sequences at the branches of a curve singularity coincide with those of its Arf closure;
moreover, by using Arf rings, it is possible to characterize all admissible multiplicity trees of curve singularities (see \cite{BDF0}).
Arf local rings are rings with minimal multiplicity, since they can be defined as rings for which every regular integrally closed ideal is stable. For this class of rings, we provide a characterization in terms of reflexive fractional  ideals, showing that these are exactly, up to isomorphism, the integrally closed ideals lying between $R:\overline{R}$ and $R$
(see Theorem \ref{Characterization Arf}). Note that one implication of this theorem was proved in a different manner in \cite{IK}, while in our paper it is a straightforward consequence of Lemma \ref{I:I-2}.
We also recall that a different characterization of Arf rings involving reflexive modules is given in \cite{D}.

\section{Preliminaries}

Let $(R, \mathfrak m)$ be a local, noetherian, one-dimensional, Cohen-Macaulay ring. Denote by $Q(R)$ its total ring of fractions and by $\overline R$ the integral closure of $R$ in $Q(R)$. Throughout the article, we assume that the residue field $R/\mathfrak m$ is infinite and that $R$ is analytically unramified, i.e., $\hat R$ is reduced.
Under our assumptions, the latter is equivalent to $\overline R$ being a finite $R$-module, see \cite[Theorem 10.2]{M}.

By a (fractional) ideal, we always mean a regular (fractional) ideal, that is, one containing a regular element. We recall that two fractional ideals $I,J$ are isomorphic if and only if there exists an invertible element $x \in Q(R)$
such that $I=xJ$. When we write $I:J$ we always mean that the colon is taken in $Q(R)$.

From our assumptions, it follows that $R$ is reduced, since its completion is reduced, and that the conductor $C=(R:\overline R) \neq (0)$ is a regular ideal of $R$; moreover, being $\hat R$ reduced, so generically Gorenstein, there exists a canonical ideal $\omega_R$, i.e. a fractional ideal such that $\omega_R:(\omega_R:I)=I$ for any fractional ideal $I$, see \cite[Satz 6.21]{HK}. 

Recall that $J$ is said to be a {\it reduction} of $I \subseteq R$, if $J\subseteq I$ and $JI^n=I^{n+1}$ for some $n \in \mathbb{N}$. Moreover, $J$ is said to be {\it minimal} if there are no other reductions of $I$ properly contained in $J$. It is well known that, under our hypotheses, any minimal reduction of an $\m$-primary ideal $I$ is
principal.

\begin{remark}\label{m}\rm
		If $(x)$ is a minimal reduction of an ideal $I$, we have that $\overline I=\overline{(x)}=x\overline R \cap R$. In fact, denoting by $P_1, \dots , P_h$ the minimal primes of $R$, $R  \hookrightarrow R/P_1 \times \dots \times  R/P_h$, $Q(R) \cong Q(R/P_1) \times \dots \times Q(R/P_h)$ and $\overline R \cong \overline{(R/P_1)} \times \dots \times \overline{(R/P_h)}$, where $\overline{(R/P_i)}$ is the integral closure of $R/P_i$ in its quotient field $Q(R/P_i)$ (see \cite[Proposition 5.17]{G}). For any $M$ maximal ideal of $\overline R$, the localization $\overline R_M$ is a DVR, hence $(x\overline R)_M=x\overline R_M$ is an integrally closed ideal, and, by \cite[Proposition 1.1.4]{HS}, it follows that $x\overline R$ is also integrally closed.
		Therefore, by \cite[Proposition 1.6.1]{HS}, the integral closure of $(x)$ as ideal of $R$ equals $x\overline R \cap R$, and, being $(x)$ a reduction of $I$, we get $I \subseteq \overline{(x)}=x\overline R \cap R$; finally, by $x \in I$, we also get the equality $\overline I=\overline{(x)}$.
\end{remark}

By the previous remark, it follows that the ideal $(x)$ is a minimal reduction of an ideal $I$ if and only if $x^{-1}I \subseteq \overline R$, and in this case we also say that $x$ is a minimal reduction of $I$.
        
For any fractional ideal $I$, let $y\in R $ be such that $yI \subseteq R$, and let $(x)$ be  a minimal reduction of $yI$.
We will call $x/y$ a \emph{reduction} of the fractional ideal $I$. Note that $x/y \in I$.
The choice of $x$, $y$ implies that
 $$R \subseteq \left(\frac{y}{x}\right)I \subseteq \overline R.$$
 In particular, up to isomorphism, we can assume that the canonical ideal $\omega_R$ is such that $R\subseteq \omega_R \subseteq \overline R$, and we will maintain this assumption throughout the paper.

In the following remark we list some properties that we will use in sequel.
 
\begin{remark} \label{Equalities} \rm  The following equalities hold:
	\begin{enumerate}
		\item $\omega_R:\omega_R=R$;
		\item $\omega_R : \overline{R}=C$;
		\item $\omega_R:C=\overline{R}$;
		\item $R:C=\overline R$.
	\end{enumerate}
The first equality holds for any canonical ideal, since $R=\omega_R:(\omega_R:R)=\omega_R:\omega_R$. 
As for the second one, we have $1\in \omega_R \subseteq \overline{R}$, so that $\overline R=\omega_R\overline R$ and, therefore, $\omega_R : \overline{R}=\omega_R : \omega_R\overline{R}=(\omega_R:\omega_R):\overline R= R:\overline R=C$.
Moreover, dualizing item (2) by $\omega_R$, we immediately obtain the equality (3).

To check the last equality, we dualize by $\omega_R$, write $R=\omega_R:\omega_R$ and use item (2). Hence, the equality is equivalent to 
$\omega_R:((\omega_R:\omega_R):C)=\omega_R:\overline R$.
On the other hand, 
$\omega_R:((\omega_R:\omega_R):C)=\omega_R:(\omega_R:\omega_R C)=\omega_R C=C=\omega_R:\overline R$, so the thesis is verified.
\end{remark}

We conclude this section by showing a result that will allow us to focus only on the fractional ideals between $C$ and $\omega_R$.

\begin{lemma}\label{FractIdeals} Let $I$ be a fractional ideal and let $z$ be a reduction of $\omega_R:I$.
Then $$C \subseteq zI \subseteq \omega_R.$$
In particular, every fractional ideal is isomorphic to a fractional ideal lying between $C$ and $\omega_R$.
\end{lemma}

\begin{proof}
    Since $z$ is a reduction of $\omega_R:I$, we have $z^{-1}(\omega_R:I) \subseteq \overline{R}$. Then, dualizing and using (2) of Remark \ref{Equalities}, we get 
    \[
    C=\omega_R:\overline{R}\subseteq \omega_R:(z^{-1}(\omega_R:I))=zI.
    \]
    The other containment holds because $z \in (\omega_R:I)$.
\end{proof}

\section{Divisorial closure  versus integral closure}

A fractional ideal $I$ is said to be \emph{reflexive} (or \emph{divisorial}) if $R:(R:I)=I$.
We gather some well-known results that we will frequently use without further reference.

\begin{remark}\label{facts}\rm
\begin{enumerate}
\item It follows immediately from the definition that a principal ideal is reflexive. 
\item  Every integrally closed ideal is reflexive. This was first proved in \cite[Proposition 2.14]{CHKV} when $R$ is a domain, while in general is a consequence of \cite[Proposition 3.9]{DMS}.
\item In particular, the maximal ideal of a local ring is integrally closed, so $\m$ is always reflexive. 
\item Reflexivity is preserved under isomorphism. In fact, if $I$ and $J$ are two isomorphic fractional ideals, then $J=yI$ for some $y \in Q(R)$ invertible in $Q(R)$. Moreover, $R:(R:yI)=R:y^{-1}(R:I)=y(R:(R:I))$, and thus $I$ is reflexive if and only if $J=yI$ is reflexive.
\end{enumerate}
\end{remark}

The next result shows that every reflexive fractional ideal is isomorphic to an ideal containing $C$. This is a more precise version of \cite[Theorem 3.5]{DMS}.

\begin{proposition} \label{Every ideal is isomorphic to a proper ideal} 
Let $I$ be a fractional reflexive ideal and let $z$ be a reduction of $\omega_R:I$. Then, there exists an invertible element $w \in \overline R$ such that  $C \subseteq wzI \subseteq R$. In particular, $I$ is isomorphic to an ideal containing $C$.
\end{proposition}

\begin{proof}
Let $w$ be a reduction of $R:zI$. It follows immediately that $wzI \subseteq R$. Moreover, $w \in R:zI \subseteq \omega_R:zI = z^{-1}(\omega_R:I) \subseteq \overline R$, where the last containment holds since $z$ is a reduction of $\omega_R:I$.

We now claim that $w$ is invertible in $\overline R$.
Assume the contrary. It follows that $w^{-1}\in Q(R)\setminus \overline R$. Moreover, by $zI \subseteq \omega_R$, we get $zI:C\subseteq \omega_R:C=\overline R$ (see Remark \ref{Equalities}). Therefore, $w^{-1} \notin zI:C$ and thus there exists $c \in C$ such that $w^{-1}c \notin zI$. 

Let $b$ be any element of $R:zI$; we have
$(w^{-1}c)b=c(w^{-1}b)  \in c\overline R$, since $w$ is a reduction of $R:zI$. Moreover, remembering that $C$ is an ideal of $\overline R$, 
we get $c\overline R \subseteq C$.
Thus, keeping in mind that $I$ is reflexive, we obtain $
w^{-1}c \in C:(R:zI) \subseteq R:(R:zI)=z(R:(R:I))=zI$, which is a contradiction.
Therefore, $w$ is invertible in $\overline R$. 

Finally, by Lemma \ref{FractIdeals}, we have $zI \supseteq C$, which implies $wzI \supseteq wC=C$, where the last equality holds because $w$ is a unit in $\overline R$ and $C$ is an ideal of $\overline{R}$.
\end{proof}

In light of the previous result, we now restrict our attention to ideals containing the conductor.

\begin{lemma}\label{I:I-2}
	Let  $I$ be a proper ideal containing $C$ and $J$ be an integrally closed ideal such that $I \subseteq J$. Then: 
	\begin{enumerate}
		\item $R:I=J:I$.
		\item If $I$ is non-principal and reflexive, $J:J \subseteq I:I$.
        \item If $I$ is non-principal and reflexive, and $(x)$ is a minimal reduction of $I$, then $x(\overline{I}:\overline{I}) \subseteq I$. 
	\end{enumerate}
\end{lemma}

\begin{proof}
	(1) Let $(x)$ be a minimal reduction of $J$; hence, $J=x\overline R \cap R$ by Remark \ref{m}.
	We always have $R:I \supseteq J:I$. Conversely, let $y \in R:I$, which means $yi \in R$ for any $i \in I$. Since $I \subseteq J$, we get $ix^{-1} \in \overline R$, whereas $y \in R:I \subseteq R: C= \overline R$ implies that $yix^{-1} \in \overline R$; thus, $yi \in x\overline R$. Therefore, $yi \in x\overline R \cap R=J$ and the thesis follows immediately.
    
    \noindent (2) Let $y \in (J:J)$ and $i \in I$. Since $R:I=J:I$ by point (1), for every $\lambda \in R:I$ we have $yi\lambda \in y J \subseteq J \subseteq R$ and, then, $yi \in R:(R:I)=I$ implies $y \in (I:I)$. \\
    \noindent (3) Using (2), we obtain that $\overline{I}:\overline{I} \subseteq I:I$ and, since $x \in I$, it immediately follows that $x(\overline{I}:\overline{I})\subseteq x(I:I)\subseteq I$.  
\end{proof}

\begin{example} \rm \label{counterexample I:I-2}
Under the assumptions of Lemma \ref{I:I-2}, the inclusion $J:J \subseteq I:I$ does not imply that $I$ is reflexive. Indeed, consider the semigroup ring $R=\mathbb{K}[[t^7,t^{10},t^{13}]]$, where $\mathbb{K}$ is a field and $t$ is an indeterminate. Let $I=(t^7,t^{13})$ and let $J$ be the maximal ideal of $R$. Since both $J:J$ and $I:I$ are monomial fractional ideals, it is not difficult to verify that they coincide; in fact, they are equal to the fractional ideal generated by $1,t^{29},t^{32}$. However, $I$ is not reflexive since $t^{10} \in (R:(R:I)) \setminus I$. 

Moreover, this example also shows that the converse of Lemma \ref{I:I-2} (3) does not hold. Indeed, note that the integral closure of $R$ is $\overline{R}=\mathbb{K}[[t]]$, and $x=t^7$ is a minimal reduction of $I$. Therefore, $\overline{I}=x\overline{R}\cap R=(t^7,t^{10},t^{13})$ is the maximal ideal of $R$. Since $\overline{I}:\overline{I}$ is generated by $1,t^{29},t^{32}$, it follows that $x(\overline{I}:\overline{I}) \subseteq I$. 
\end{example}

The operation that associates to an ideal $I$ the ideal $R:(R:I)$ is the classical divisorial closure, also known as $v$-closure.
In the following proposition, we compare divisorial and integral closures.

\begin{proposition}\label{main reflexivity}
Let $(x)$ be a minimal reduction of an ideal $I$. It holds that $R:(R:I) \subseteq \overline I$. Moreover, if $C \subseteq I$ and $I$ is not principal, then $x(\overline{I}:\overline{I}) \subseteq R:(R:I) \subseteq \overline I$.
\end{proposition}

\begin{proof}
Keeping in mind that $\omega_R I\subseteq \omega_R$, we have $R:(R:I)= R:((\omega_R:\omega_R):I)=R:(\omega_R:\omega_RI)\subseteq R:(\omega_R:\omega_R)=R:R=R$. Therefore, $R:(R:I)$ is an ideal of $R$, and we claim that $(x)$ is a minimal reduction of it.

To this aim, we first notice that $x \in I \subseteq R:(R:I)$.
Moreover, since $(x)$ is a minimal reduction of $I$, it holds $x^{-1} I \subseteq \overline{R}$. In particular, $x^{-1} I \omega_R \subseteq \overline{R}$ and then by Remark \ref{Equalities} we obtain
\[
C=\omega_R : \overline{R} \subseteq \omega_R : (x^{-1} I \omega_R)=(\omega_R:\omega_R):x^{-1}I=R:x^{-1}I.
\]
Thus, $C \subseteq  R:x^{-1}I$ implies $x^{-1}(R:(R:I))=R:(R:x^{-1}I) \subseteq R:C=\overline{R}$, so that $(x)$ is a minimal reduction of $R:(R:I)$.

As a consequence, $R:(R:I) \subseteq \overline{(x)}=\overline{I}$. Moreover, if $C \subseteq I$, we have $C \subseteq I \subseteq R:(R:I)\subseteq R$ and $R:(R:I)$ is clearly reflexive; then, the thesis follows immediately by the point (3) of Lemma \ref{I:I-2}.
\end{proof}

Note that the inclusion $R:(R:I) \subseteq \overline{I}$ was already proved in \cite[Proposition 2.14]{CHKV} when $R$ is a domain.

\begin{example} \rm
In the second part of the previous proposition, we cannot remove the assumption that $C \subseteq I$. For instance,
let $R=\mathbb{K}[[t^3,t^5,t^7]]$ and $I=t^2 \m=(t^5,t^7,t^9)$. Since $I$ is isomorphic to $\m$, it is reflexive. Moreover, being an element of $I$ of minimal valuation, $(t^5)$ is a minimal reduction of $I$ and $\overline{I}=C=t^5\mathbb{K}[[t]]$. Therefore, $t^5(C:C)=C \not\subseteq I=R:(R:I)$.    
\end{example}

\section{Reflexive ideals and almost Gorenstein rings}

We begin this section with a remark that provides a sufficient condition for an ideal to be reflexive.

\begin{remark} \rm \label{omega in I:I}
If $I$ is a fractional ideal such that $I:I$ contains $\omega_R$, then $I$ is reflexive. Indeed,
\begin{gather*}
R:(R:I)=(\omega_R:\omega_R):((\omega_R:\omega_R):I)=\omega_R:(\omega_R(\omega_R:\omega_R I))= \\
=(\omega_R:(\omega_R:\omega_R I)):\omega_R=\omega_R I:\omega_R \subseteq (I:I)I:\omega_R \subseteq I:\omega_R \subseteq I.
\end{gather*}
\end{remark}

In general, the converse of the previous remark does not hold, as the next example shows.

\begin{example} \rm \label{Jager}
Consider the numerical semigroup ring $R=\mathbb{K}[[t^5, t^7,t^9]]$, where $\mathbb{K}$ is a field, and its ideal $I=(t^7,t^9,t^{10})$. It is easy to see that $I=t^7\mathbb K[[t]]\cap R$, hence it is integrally closed and then reflexive. In order to show that $\omega_R \nsubseteq I:I$, we need to compute $\omega_R$. Fortunately, for a numerical semigroup ring there is an easy description of $\omega_R$. Let $S=\langle 5,7,9 \rangle=\{0,5,7,9,10,12,14 \rightarrow\}$ be the associated semigroup of $R$, i.e. the semigroup consisting of the exponents of the monomials in $R$, where $\rightarrow$ means that all the integers greater than $14$ are in $S$. The integer $\mathrm{F}(S)=\max(\mathbb{N} \setminus S)=13$ is called Frobenius number of $S$. Let $K(S)=\{\mathrm{F}(S)-s \mid s \notin S\}=\{0,2,5,7,9,10,11,12,14 \rightarrow\}$. By \cite[Satz 5]{J}, the fractional ideal $\omega_R=\mathbb{K}[[t^s \mid s \in K(S)]]$ is a canonical ideal of $R$ such that $R \subseteq \omega_R \subseteq \overline{R}=\mathbb{K}[[t]]$.
Therefore, $t^2 \in \omega_R$ but $t^2 t^9 \notin I$, and thus $\omega_R \nsubseteq I:I$.  
\end{example}

The first goal of this section is to show that the converse of Remark \ref{omega in I:I} holds for every non-principal (fractional) ideal exactly when $R$ is an almost Gorenstein ring. We recall that, with our assumptions, $R$ is said to be {\it almost Gorenstein} when $\m\omega_R \subseteq \m$.

\begin{remark} \label{I ideal of B} \rm
If $I$ is a non-principal reflexive fractional ideal of $R$, then $\m:\m \subseteq I:I$. Indeed, by Proposition \ref{Every ideal is isomorphic to a proper ideal}, there exists $\alpha \in Q(R)$ such that $C \subseteq \alpha I \subseteq R$. Moreover, since $I$ is not principal, $\alpha I \subseteq \m$. Therefore, Lemma \ref{I:I-2} implies that $\m:\m \subseteq \alpha I : \alpha I=I:I$. 
\end{remark}

Now we can provide a new characterization of the almost Gorenstein property of $R$.

\begin{proposition} \label{almost Gorenstein}
The following statements are equivalent:
\begin{enumerate}
\item $R$ is almost Gorenstein.
\item A non-principal fractional ideal $I$ of $R$ is reflexive if and only if $\omega_R \subseteq I:I$.
\item A non-principal ideal $I$ of $R$ is reflexive if and only if $\omega_R \subseteq I:I$.
\end{enumerate}
\end{proposition}

\begin{proof}
(1) $\Rightarrow$ (2) If $\omega_R \subseteq I:I$, we have already seen in Remark \ref{omega in I:I} that $I$ is reflexive. Conversely, by the definition of almost Gorenstein ring and by the previous remark, we obtain $\omega_R \subseteq \m:\m \subseteq I:I$.  \\
(2) $\Rightarrow$ (3) Clear. \\
(3) $\Rightarrow$ (1) Since $\m$ is reflexive, by assumption we have $\omega_R \subseteq \m:\m$, which means that $R$ is almost Gorenstein.
\end{proof}

The next goal of the section is to provide another characterization of almost Gorenstein rings $R$, which involves reflexive ideals of $R$ and ideals of the overring $\m:\m$.

Let $I$ and $H$ be two fractional ideals of $R$. 
We always have the inclusion $I \subseteq H:(H:I)$, and we say that $I$ is $H$-reflexive when the equality holds; so, when we simply write reflexive, we are omitting $R$.

\begin{proposition} \label{H-reflexivity}
If $I$ is a proper ideal and $J$ is an integrally closed ideal such that $C \subseteq I \subseteq J$, then $R:(R:I)=J:(J:I)$.
In particular, $I$ is reflexive if and only if $I$ is $J$-reflexive.
\end{proposition}

\begin{proof}
By Lemma \ref{I:I-2}, we have $R:I=J:I$, so we obtain $J:(J:I)=J:(R:I) \subseteq R:(R:I)$.
	
	Conversely, the inclusion $R:(R:I) \subseteq J:(J:I)$ is equivalent to
	$(R:I)(R:(R:I))\subseteq J $. It is obviously true that $(R:I)(R:(R:I))\subseteq R$. Hence, given a minimal reduction $(x)$ of $J$, if we show that 
    \begin{equation} \label{equat}
    x^{-1}(R:I)(R:(R:I)) \subseteq \overline R,
    \end{equation}
    we obtain $(R:I)(R:(R:I)) \subseteq x \overline R \cap R=J$ because $J$ is integrally closed. 
    On the other hand, $R:I \subseteq R:C  = \overline R$, while $x^{-1}(R:(R:I))\subseteq x^{-1}(R:(R:J))=x^{-1}J\subseteq \overline R$ because $J$ is integrally closed and then reflexive; hence, (\ref{equat}) holds.
\end{proof}

We recall that $\m:\m$ is a ring. In general, a fractional ideal $I$ of $R$ is not a fractional ideal of $\m:\m$, but this property and reflexivity can be used to characterize the almost Gorenstein property of $R$.

\begin{theorem} \label{almost Gorenstein 2}
Let $B=\m:\m$. The following assertions are equivalent:
\begin{enumerate}
    \item $R$ is almost Gorenstein.
    \item A non-principal fractional ideal $I$ of $R$ is reflexive if and only if $I$ is a fractional ideal of $B$.
    \item A non-principal ideal $I$ of $R$ is reflexive if and only if $I$ is an ideal of $B$.
\end{enumerate}
\end{theorem}

\begin{proof}
We assume that $R$ is not a DVR, otherwise $R= B$ and there is nothing to prove. It is well known that $R \neq B$ in this case. Indeed, when $R$ is not a DVR, it always holds $R:\m=\m:\m$ since if there exists $ t \in (R:\m) \setminus (\m:\m)$, then $t\m=R$ implies that $R$ is a DVR because $\m \cong R$. Therefore, if $R=B$, one has $R=\m:\m=R:\m$ and then $R=R:R=R:(R:\m)=\m$ yields a contradiction.

(1) $\Rightarrow$ (2) Note that $I$ is a fractional ideal of $B$ if and only if $\m:\m \subseteq I:I$. 
If $I$ is reflexive, then Remark \ref{I ideal of B} immediately implies that $I$ is a fractional ideal of $B$. Conversely, if $I$ is an ideal of $B$, we have $\omega_R \subseteq \m:\m \subseteq I:I$, and then $I$ is reflexive by Proposition \ref{almost Gorenstein}. \\
(2) $\Rightarrow$ (3) Clear. \\
(3) $\Rightarrow$ (1) By \cite[Proposition 4.2]{DS}, $R$ is almost Gorenstein if and only if $\m$ is a canonical ideal of $B$, which means that every fractional ideal of $B$ is $\m$-reflexive. So, let $I$ be a fractional ideal of $B$. Clearly, $I$ is also a fractional ideal of $R$, and then there exists a unit $\alpha$ of $Q(R)$ such that $\alpha I \subseteq R$.
If $\alpha I$ is principal, then there exists $\lambda \in Q(R)$ such that  $\lambda \alpha I =R$, and we obtain $B \subseteq I:I=\lambda \alpha I:\lambda \alpha I=R:R=R$, which is not possible because $R$ is not a DVR; thus, $\alpha I$ is not principal.

 Moreover, since $I$ is an ideal of $B$, so is $\alpha I$. Therefore, by hypothesis, $\alpha I$ is reflexive. By Proposition \ref{Every ideal is isomorphic to a proper ideal}, there exists a unit $\lambda$ of $Q(R)$ such that $C \subseteq \lambda \alpha I \subseteq \m$, since $\lambda \alpha I$ is not principal and then different from $R$. 
 By Proposition \ref{H-reflexivity}, $\lambda \alpha I$ is $\m$-reflexive, and this easily implies that also $I$ is $\m$-reflexive because $\lambda \alpha$ is a unit of $Q(R)$.
\end{proof}

\section{Reflexive ideals and rings with minimal multiplicity}

Let $\nu(R)$ and $e(R)$ denote the embedding dimension and the multiplicity of $R$ respectively. Under our assumptions, we have $\nu(R) \leq e(R)$, and $R$ is said to have {\it minimal multiplicity} if equality holds. Moreover, $R$ has minimal multiplicity if and only if its maximal ideal satisfies the equality $x \m=\m^2$, for some $x \in \m$ (see \cite[Corollary 1.10 and Lemma 1.11]{L}). More generally, an ideal $I$ for which $x I=I^2$, for some $x \in I$, is called {\it stable}; thus, $R$ has minimal multiplicity if and only if $\m$ is stable. Notice that the stability of an ideal is invariant under isomorphisms, and, if $(x)$ is a reduction of $I$, the ideal $I$ is stable if and only if $x(I:I)=I$; see \cite{L} for more details about stable ideals.

We begin this section, providing a characterization of rings with minimal multiplicity by means of reflexivity. We start by establishing some preliminary results.

\begin{remark}\label{transistivity}\rm
     Reflexivity is a transitive property; more precisely, if $I,H,F$ are fractional ideals, $I$ is $H$-reflexive and $H$ is $F$-reflexive, then $I$ is also $F$-reflexive. In fact, we have
$$
\begin{aligned}
    F:(F:I)&=F:(F:(H:(H:I)))=F:(F:((F:(F:H)):(H:I)))= \\
    &= F:(F:(F:(F:H) (H:I)))
=F:((F:H)(H:I))= \\ &= (F:(F:H)):(H:I)=H:(H:I)=I.
\end{aligned}
$$
\end{remark}

\begin{remark} \label{I:I} \rm
If $I$ is an integrally closed ideal containing $C$, then $I:I$ is reflexive. In fact, $I:I=R:I$ by Lemma \ref{I:I-2}, and $R:(R:(I:I))=R:(R:(R:I))=R:I=I:I$. 
\end{remark}

\begin{lemma}\label{1 in X}
Let $\m$ be the maximal ideal of $R$ and set $B=(\m:\m)$. Then, the set $\m(\m:\m^2)$ is an ideal of $B$, and the maximal ideal $\m$ of $R$ is stable if and only if $\m(\m:\m^2)=B$.
\end{lemma}

\begin{proof}
  If $\m$ is stable and $(x)$ is a minimal reduction of $\m$, then $\m^2=x\m$ or, equivalently, $B=\m:\m=x^{-1}\m$. Therefore, $\m(\m:\m^2)=x^{-1}\m(\m:\m)=B\cdot B=B$, where the last equality holds since $B$ is a ring.

Conversely, let $r$ be the reduction number of $\m$, that is $x\m^r =\m^{r+1}$ but $x\m^{r-1} \subsetneq \m^r$, and assume by contradiction that $r>1$.
Since $1 \in \m(\m:\m^2)$, it follows that 
$$\m^r \subseteq\m^{r+1}(\m:\m^2)=x\m^r(\m:\m^2)=x\m^{r-2}\m^2(\m:\m^2) \subseteq 
x\m^{r-2}\m = x\m^{r-1},$$ 
which yields a contradiction.
\end{proof}

In the following theorem, we use reflexivity to characterize rings with minimal multiplicity in a manner similar to Theorem \ref{almost Gorenstein 2}.

\begin{theorem} \label{minimal multiplicity}
	Let $B=(\m:\m)$.
    The following assertions are equivalent:
    \begin{enumerate}
        \item $R$ has minimal multiplicity.
        \item  For any non-principal fractional ideal $I$ of $R$, $I$ is reflexive if and only if $I$ is a $B$-reflexive fractional ideal of $B$. 
        \item $\m$ is $B$-reflexive.
    \end{enumerate}
\end{theorem}

\begin{proof}
(1) $\Rightarrow$ (2)
By Remark \ref{I:I}, we know that $B$ is reflexive as $R$-fractional ideal. Hence, by Remark \ref{transistivity}, every $B$-reflexive fractional ideal of $B$ is reflexive.

Conversely, if $I$ is a non-principal reflexive fractional ideal of $R$, then $\m:\m \subseteq I:I$ by Remark \ref{I ideal of B},
which means that $I$ is a fractional ideal of $B$. Moreover, since $\m$ is an integrally closed ideal, $I$ is $\m$-reflexive, by Proposition \ref{H-reflexivity}. Therefore, as reflexivity is a transitive property (cf. Remark \ref{transistivity}), it is enough to show that $\m$ is $B$-reflexive. 
Thus, we have to compute $$B:(B:\m)=(\m:\m):((\m:\m):\m)=(\m:\m):(\m:\m^2)=\m:(\m(\m:\m^2))$$ and to show that this equals $\m$.
By hypothesis, $\m$ is stable, and then $\m(\m:\m^2)=B$ by Lemma \ref{1 in X}.
Hence, $B:(B:\m)=\m:B=\m$, since $\m$ is an ideal of $B$.

\noindent (2) $\Rightarrow$ (3) If $\m$ is principal, then $R$ is a DVR and the statement is trivial because $R=B$. Otherwise, it is enough to apply the hypothesis. 

\noindent (3) $\Rightarrow$ (1)
 Assume by contradiction that $\m$ is not stable; by Lemma \ref{1 in X}, this is equivalent to $\m(\m:\m^2) \subsetneq B$, which implies that $\m(\m:\m^2)$ is contained in a maximal ideal $\mathfrak n$ of $B$. 

We first notice that $C\subsetneq C:\n$. In fact, it is enough to show that this holds when we localize at $\n$, then we may assume that $B$ is local. In this case, since $C$ is regular and $\n$ is the unique prime ideal with positive height, $\n$ is the radical of $C$. Therefore, there exists $a \in \mathbb{N}$ such that $\n^a \subseteq C$ and we assume that $a$ is the minimum for which this happens. Hence, $\n^{a-1} \subseteq C:\n$ but $\n^{a-1} \not\subseteq C$, which implies $C\subsetneq C:\n$.

We also notice that $\m:\overline R=C$; indeed, $\m:\overline{R}\subseteq R:\overline R=C$ and, conversely,   $C\overline R=C\subseteq \m$.
Moreover, as computed in the reverse implication, we have $B:(B:\m)=\m:(\m(\m:\m^2))$. 
Therefore, we obtain
 $$(B:(B:\m)):\overline R=(\m:(\m(\m:\m^2))):\overline R \supseteq (\m :\mathfrak n):\overline R=
 (\m:\overline R):\mathfrak n=C:\mathfrak n \supsetneq C = \m:\overline R$$

Hence, $\m:\overline R \subsetneq (B:(B:\m)):\overline R$, but, this leads to a contradiction because $\m=(B:(B:\m))$ by assumption.
\end{proof}

The following examples clarify the difference between the characterizations given in Theorem \ref{almost Gorenstein 2} and Theorem \ref{minimal multiplicity}.

\begin{example} \rm Let $\mathbb{K}$ be a field and $t$ be an indeterminate. \\
(1) Let $R=\mathbb{K}[[t^4,t^5,t^7]]$, and $\m=(t^4,t^5,t^7)$. As in Example \ref{Jager}, one can see that the fractional ideal $\omega_R$ generated by $1$ and $t^3$ is a canonical module of $R$, and $\m \omega_R =\m$; then, $R$ is almost Gorenstein. 
The ideal $I=(t^5,t^7,t^8)$ of $R$  is reflexive. By Theorem \ref{almost Gorenstein 2}, $I$ is also an ideal of $B=\m:\m=\mathbb{K}[[t^3,t^4,t^5]]$, but it is not reflexive as ideal of $B$ because $t^6 \in (B:(B:I)) \setminus I$. This means that $R$ has not minimal multiplicity by Theorem \ref{minimal multiplicity}. \\
(2) Let $R=\mathbb{K}[[t^3,t^7,t^8]]$, $\m=(t^3,t^7,t^8)$ and $I=(t^6,t^8,t^{10})$. It is straightforward to see that $I$ is also an ideal of $B=\m:\m=\mathbb{K}[[t^3,t^4,t^5]]$, but $I$ is not reflexive because $t^7\in (R:(R:I)) \setminus I$. By Theorem \ref{almost Gorenstein 2}, $R$ is not almost Gorenstein. On the other hand, $R$ has minimal multiplicity because $\m^2=t^3\m$. 
\end{example}

The ring $R$ is called {\it Arf} if every integrally closed ideal of $R$ is stable; in particular, an Arf ring has minimal multiplicity because $\m$ is stable. This notion was originally introduced by Lipman with a different definition, but he also proved in \cite[Theorem 2.2]{L} that the definition we use here is equivalent. In the last part of this section we will provide two characterizations of Arf rings in terms of reflexive ideals.

\begin{lemma}
    Let $I \subseteq J$ be two proper ideals of $R$ containing $C$. If $J=\alpha I$ for some regular $\alpha \in Q(R)$, then  $I=J$.
\end{lemma}

\begin{proof}
    Since $I \subseteq \alpha I=J$, we get, inductively $\alpha^n I \subseteq \alpha^{n+1} I$ for any $n \geq 0$. Moreover, notice that $\alpha \in \overline R$, since $\alpha C \subseteq \alpha I \subset R$ and so $\alpha \in R:C=\overline R$.
    Thus, we have the following increasing chain of $R$-submodules of $\overline R$, which is a finite (and hence noetherian) $R$-module,
    $$
    I \subseteq \alpha I \subseteq \alpha^2 I \subseteq \dots \subseteq \alpha^n I\subseteq \dots
    $$
    that is stationary. Therefore, 
    $\alpha^h I=\alpha^{h+1}I$, for some $h \geq 0$, and, being $\alpha$ regular, we obtain $I=\alpha I=J$.
\end{proof}

In \cite[Theorem 3.4]{IK}, it is proved that in an Arf ring every reflexive ideal is isomorphic to an integrally closed ideal. The implication (1) $\Rightarrow$ (3) in the next result is a more precise version of this fact, and, in this case, we also show that the converse holds.

\begin{theorem} \label{Characterization Arf}
The following assertions are equivalent:
\begin{enumerate}
    \item The ring $R$ is Arf.
    \item Every reflexive ideal $I$ such that $C \subseteq I \subseteq R$ is integrally closed.
    \item Every reflexive fractional ideal is isomorphic to an integrally closed ideal $I$ with $C \subseteq I \subseteq R$.
\end{enumerate}
\end{theorem}

\begin{proof}
(1) $\Rightarrow$ (2) Let $I$ be a reflexive ideal such that $C \subseteq I \subseteq R$. Given a minimal reduction $(x)$ of $I$, we get $x(\overline{I}:\overline{I}) \subseteq I \subseteq \overline{I}$ by Lemma \ref{I:I-2}. On the other hand, $\overline{I}$ is integrally closed and $R$ is Arf, then $\overline{I}$ is stable, i.e., $x(\overline{I}:\overline{I})=\overline{I}$; this implies that $I=\overline{I}$, and then $I$ is integrally closed.  

(2) $\Rightarrow$ (1) We need to prove that every integrally closed fractional ideal $J$ is stable. Since $J$ is integrally closed, it is reflexive, and Proposition \ref{Every ideal is isomorphic to a proper ideal} asserts that $J$ is isomorphic to a reflexive ideal $I$ such that $C \subseteq I \subseteq R$, which is integrally closed by assumption. It is enough to show that $I$ is stable.

Since $I$ is integrally closed, it follows that $I = x \overline{R} \cap R = \overline{(x)}$ for a minimal reduction $(x)$ of $I$. We know that $x^{-1} I \subseteq \overline{R}$ and that $C$ is an ideal of $\overline{R}$, thus $x^{-1} I C \subseteq C \subseteq I$ and this implies $C \subseteq I: x^{-1} I = x(I:I) \subseteq I \subseteq R$. Moreover, $x(I:I)$ is reflexive because $(I:I)$ is reflexive by Remark \ref{I:I}. It follows that $x(I:I)$ is integrally closed by assumption. 
Furthermore, we have $I:I=R:I$ by Lemma \ref{I:I-2}, and \cite[Lemma 2.3]{DS} implies that $(1)$ is a reduction of $I:I$ because $R:I \subseteq R:C=\overline R$; then, $(x)$ is a minimal reduction of $x(I:I)$. Hence, $\overline{(x)}=x(I:I) \subseteq I = \overline{(x)}$ implies that $x(I:I)=I$, i.e., $I$ is stable.

(2) $\Rightarrow$ (3)  It is enough to apply Proposition \ref{Every ideal is isomorphic to a proper ideal} and use the assumption. 

(3) $\Rightarrow$ (2) Let $I$ be a reflexive ideal of $R$ with $C \subseteq I \subseteq R$. By assumption, $I$ is isomorphic to an integrally closed ideal $J$ with $C \subseteq J \subseteq R$, i.e., there exists $\alpha \in Q(R)$ such that $J=\alpha I$ is integrally closed and $C \subseteq \alpha I \subseteq R$.

Since $\alpha I \subseteq R$ and $C \subseteq I$, it follows that $\alpha C \subseteq R$, i.e., $\alpha \in R:C = \overline{R}$, where the last equality holds by Remark \ref{Equalities}; thus, $\alpha$ is in $\overline{R}$. Moreover, $C \subseteq \alpha I$ implies that $\alpha^{-1} C \subseteq R$, i,e., $\alpha^{-1} \in (R:C)=\overline{R}$. 
Therefore, $\alpha$ is a unit in $\overline{R}$. 

Let $(x)$ be a minimal reduction of $I$; notice that $\alpha x$ is a minimal reduction of $\alpha I$, because $x I^n=I^{n+1}$ implies $\alpha x (\alpha I)^n=(\alpha I)^{n+1}$. Thus, $\alpha I=\overline{\alpha I}=\alpha x \overline{R} \cap R= x \overline{R} \cap R= \overline{I}$, where the third equality holds because $\alpha$ is a unit in $\overline{R}$. This means that $I \subseteq \overline{I}=\alpha I$. 
By the previous lemma we get $I=\overline I$.
\end{proof}

In the previous theorem, the condition that $I$ contains the conductor cannot be removed, as the following example shows.

\begin{example} \rm
Consider $R=\mathbb{K}[[t^3,t^7,t^8]]$. Since $R$ is a numerical semigroup ring, one can verify that it is an Arf ring by applying \cite[Theorem II.2.13 and Theorem I.3.4]{BDF}. Observe that the ideal $I=(t^6,t^{10},t^{11})$ is reflexive and not integrally closed . This is due to the fact that $I$ does not contain the conductor. However, according to Theorem \ref{Characterization Arf}, $I$ is isomorphic to $t^{-3} I$, which coincides with the maximal ideal of $R$, and is therefore integrally closed.
\end{example}

\end{document}